\documentstyle[twoside]{article}
\title{On extensions of two results due to Ramanujan}
\author{\sc Y. S. Kim$^a$,  A. K. Rathie$^b$ and R. B. Paris$^c$\\
\\
${}^a\!$ Department of Mathematics Education, Wonkwang University, Iksan, Korea\\
E-Mail: yspkim@wonkwang.ac.kr\\
${}^b\!$ Department of Mathematics, Central University of Kerala,\\ Kasaragad 671328, Kerala, India\\
E-Mail: akrathie@gmail.com\\
${}^c\!$ University of Abertay Dundee, Dundee DD1 1HG, UK\\
E-Mail: r.paris@abertay.ac.uk}
\begin{document}
\def\f#1#2{\mbox{${\textstyle \frac{#1}{#2}}$}}
\def\dfrac#1#2{\displaystyle{\frac{#1}{#2}}}
\def\boldal{\mbox{\boldmath $\alpha$}}
\newcommand{\bee}{\begin{equation}}
\newcommand{\ee}{\end{equation}}
\newcommand{\lam}{\lambda}
\newcommand{\ka}{\kappa}
\newcommand{\al}{\alpha}
\newcommand{\th}{\theta}
\newcommand{\om}{\omega}
\newcommand{\Om}{\Omega}
\newcommand{\fr}{\frac{1}{2}}
\newcommand{\fs}{\f{1}{2}}
\newcommand{\g}{\Gamma}
\newcommand{\br}{\biggr}
\newcommand{\bl}{\biggl}
\newcommand{\ra}{\rightarrow}
\newcommand{\mbint}{\frac{1}{2\pi i}\int_{c-\infty i}^{c+\infty i}}
\newcommand{\mbcint}{\frac{1}{2\pi i}\int_C}
\newcommand{\mboint}{\frac{1}{2\pi i}\int_{-\infty i}^{\infty i}}
\newcommand{\gtwid}{\raisebox{-.8ex}{\mbox{$\stackrel{\textstyle >}{\sim}$}}}
\newcommand{\ltwid}{\raisebox{-.8ex}{\mbox{$\stackrel{\textstyle <}{\sim}$}}}
\renewcommand{\topfraction}{0.9}
\renewcommand{\bottomfraction}{0.9}
\renewcommand{\textfraction}{0.05}
\newcommand{\mcol}{\multicolumn}
\date{}
\maketitle
\begin{abstract}
The aim in this note is to provide a generalization of an interesting entry in Ramanujan's Notebooks
that relate sums involving the derivatives of a function $\phi(t)$ evaluated at 0 and 1. 
The generalization obtained is derived with the help of expressions for the sum of terminating ${}_3F_2$ hypergeometric functions of argument equal to 2, recently obtained in Kim {\it et al.} [{\it Two results for the terminating ${}_3F_2(2)$ with applications}, Bull. Korean Math. Soc. {\bf 49} (2012) pp. 621--633].
Several special cases are given. In addition we generalize a summation formula to include integral parameter differences.

\vspace{0.4cm}

\vspace{0.3cm}

\end{abstract}

\vspace{0.3cm}

\begin{center}
{\bf 1. \  Introduction}
\end{center}
\setcounter{section}{1}
\setcounter{equation}{0}
\renewcommand{\theequation}{\arabic{section}.\arabic{equation}}
Two of the many interesting results stated by Ramanujan in his Notebooks are the following theorems,
which appear as Entries 8 \cite[p.~51]{B} and Entry 20 \cite[p.~36]{B}, expressing an infinite sum of derivatives of a function $\phi(t)$ at the origin to another infinite sum of its derivatives evaluated at $t=1$.
\vspace{0.25cm}

\noindent{\bf Entry 8.} \ {\it Let $\phi(t)$ be analytic for $|t-1|<R$, where $R>1$. Suppose that $a$ and $\phi(t)$ are such that the order of summation in
\[
\sum_{k=0}^\infty \frac{2^k (a)_k}{(2a)_k\,k!}\,\sum_{n=k}^\infty \frac{(-1)^n}{n!}\,(-n)_k \phi^{(n)}(1)
\]
may be inverted. Then
\bee\label{e11}
\sum_{k=0}^\infty \frac{2^k(a)_k\phi^{(k)}(0)}{(2a)_k\,k!}=\sum_{k=0}^\infty\frac{ \phi^{(2k)}(1)}{2^{2k}(a+\fs)_k\,k!}.
\ee}
\vspace{0.25cm}

\noindent{\bf Entry 20.} \ {\it Let $\phi(t)=\sum_{k=0}^\infty \phi^{(k)}(1) (t-1)^k/k!$ be analytic for $|t-1|<R$, where $R>1$. Suppose that $a$ and $b$ are complex parameters such that the order of summation in
\[
\sum_{k=0}^\infty\frac{(a)_k}{(b)_k k!}\sum_{n=k}^\infty\frac{(-1)^n}{n!} (-n)_k \phi^{(n)}(1)
\]
may be inverted. Then
\bee\label{e12}
\sum_{k=0}^\infty\frac{(a)_k\,\phi^{(k)}(0)}{(b)_k k!}=\sum_{k=0}^\infty \frac{(-1)^k (b-a)_k}{(b)_k k!} \phi^{(k)}(1).
\ee}
Berndt \cite{B} pointed out that Entry 8 can be established with the help of the results \[{}_2F_1\left[\!\!\begin{array}{c}-2n,a\\2a\end{array}\!;2\right]=\frac{(\fs)_n}{(a+\fs)_n},\qquad {}_2F_1\left[\!\!\begin{array}{c}-2n-1,a\\2a\end{array}\!;2\right]=0\]
for non-negative integer $n$, where $(a)_n=\g(a+n)/\g(a)$ denotes the Pochhammer symbol. 

In \cite{Kim}, each of the above theorems was generalized. The result in (\ref{e11}) was extended by replacing the denominatorial parameter $2a$ by $2a+j$, where $j=0, \pm 1, \dots ,\pm5$. The second result
in (\ref{e12}) was extended by the inclusion of an additional pair of numeratorial and denominatorial parameters differing by unity to produce the following theorem.
\newtheorem{theorem}{Theorem}
\begin{theorem}$\!\!\!$.\ \ Let $\phi(t)$ be analytic for $|t-1|<R$, where $R>1$. Suppose that $a$, $b$, $d$  and $\phi(t)$ are such that the order of summation in
\[\sum_{k=0}^\infty \frac{(a)_k (d+1)_k}{(b)_k (d)_k k!}\sum_{n=k}^\infty \frac{(-1)^n}{n!}(-n)_k \phi^{(n)}(1)\]
may be inverted. Then
\bee\label{e35}
\sum_{k=0}^\infty \frac{(a)_k (d+1)_k}{(b)_k (d)_k}\,\frac{\phi^{(k)}(0)}{k!}=\sum_{k=0}^\infty\frac{(-1)^k (b-a-1)_k (f+1)_k}{(b)_k (f)_k}\,\frac{\phi^{(k)}(1)}{k!},
\ee
where $f=d(b-a-1)/(d-a)$.
\end{theorem}
This result was extended to the case where a pair of numeratorial and denominatorial parameters differs by a positive integer $m$ to produce
\[
\sum_{k=0}^\infty \frac{(a)_k (d+m)_k}{(b)_k (d)_k}\,\frac{\phi^{(k)}(0)}{k!}=\sum_{k=0}^\infty\frac{(-1)^k (b-a-m)_k }{(b)_k }\,\frac{((\xi_m+1))_k}{((\xi_m))_k}\,\frac{\phi^{(k)}(1)}{k!},
\]
where the $\xi_1, \ldots , \xi_m$ are the zeros of a certain polynomial of degree $m$. 

In this note we shall similarly extend the result in (\ref{e11}) (when the parameter $2a$ is replaced by $2a+1$) by the inclusion of a pair of numeratorial and denominatorial parameters differing by unity.
For this we shall require the summations of a ${}_3F_2$ hypergeometric function of argument equal to 2
obtained\footnote{It should be noted that the right-hand side of (\ref{e18a}) is independent of the parameter $d$.} in \cite[Theorem 2]{Kim2}
\bee\label{e18a}
{}_3F_2\left[\!\!\begin{array}{c } -2n, a, d+1\\2a+1, d\end{array}\!;2\right]=\frac{(\fs)_n}{(a+\fs)_n},
\ee
and
\bee\label{e18b}
{}_3F_2\left[\!\!\begin{array}{c } -2n-1, a, d+1\\2a+1, d\end{array}\!;2\right]=\frac{(1-2a/d)}{2a+1}\,\frac{(\f{3}{2})_n}{(a+\f{3}{2})_n},
\ee
for non-negative integer $n$. Several applications are presented in Section 3.

In the final section we generalize the result given in \cite[p.~25]{B} as Entry 9:
\vspace{0.25cm}

\noindent{\bf Entry 9.} \ {\it If Re$\,(c-a)>0$ then
\bee\label{e19}
\sum_{k=1}^\infty\frac{(a)_k}{k (c)_k}=\psi(c)-\psi(c-a),
\ee
where $\psi(x)$ denotes the logarithmic derivative of\/ $\g(x)$}.
\vspace{0.25cm}

\noindent We extend this summation to include additional numeratorial and denominatorial parameters differing by positive integers. To achieve this we make use of the generalized Karlsson-Minton summation formula for a
${}_{r+2}F_{r+1}$ hypergeometric function of unit argument.

\vspace{0.6cm}

\begin{center}
{\bf 2. \ Generalization of Ramanujan's result (\ref{e11})}
\end{center}
\setcounter{section}{2}
\setcounter{equation}{0}
\renewcommand{\theequation}{\arabic{section}.\arabic{equation}}
The result to be established in this section is given by the following theorem.
\vspace{0.1cm}
\begin{theorem}$\!\!\!$.\ \ Let $\phi(t)$ be analytic for $|t-1|<R$, where $R>1$. Suppose that $a$, $d$ and $\phi(t)$ are such that the order of summation in
\[
\sum_{k=0}^\infty \frac{2^k (a)_k (d+1)_k}{(2a+1)_k (d)_k\,k!}\,\sum_{n=k}^\infty \frac{(-1)^n}{n!}\,(-n)_k \phi^{(n)}(1)
\]
may be inverted. Then
\bee\label{e21}
\sum_{k=0}^\infty \frac{2^k(a)_k (d+1)_k\phi^{(k)}(0)}{(2a+1)_k (d)_k\,k!}=\sum_{k=0}^\infty\frac{ \phi^{(2k)}(1)}{2^{2k}(a+\fs)_k\,k!}-\frac{(1-2a/d)}{2a+1} \sum_{k=0}^\infty \frac{ \phi^{(2k+1)}(1)}{2^{2k}(a+\f{3}{2})_k\,k!}.
\ee
\end{theorem}
\vspace{0.2cm}

\noindent{\it Proof.}\ \ \ Since $\phi(t)$ is analytic for $|t-1|<R$, we have
\[\phi^{(k)}(t)=\sum_{n=k}^\infty \frac{(-1)^n (-n)_k}{n!}\,\phi^{(n)}(1)\]
by suitable differentiation of the associated Taylor series.
Then
\begin{eqnarray*}
S:=\sum_{k=0}^\infty \frac{2^k(a)_k (d+1)_k}{(2a+1)_k (d)_k}\,\frac{\phi^{(k)}(0)}{k!}\!\!&=&\!\!\sum_{k=0}^\infty \frac{2^k(a)_k (d+1)_k}{(2a+1)_k (d)_k} \sum_{n=k}^\infty \frac{(-1)^n(-n)_k}{n!}\,\phi^{(n)}(1)\\
&=&\!\!\sum_{n=0}^\infty \frac{(-1)^n \phi^{(n)}(1)}{n!}\sum_{k=0}^n \frac{2^k(a)_k (d+1)_k (-n)_k}{(2a+1)_k (d)_k k!}\\
&=&\!\!\sum_{n=0}^\infty \frac{(-1)^n \phi^{(n)}(1)}{n!}\,{}_3F_2\left[\!\!\begin{array}{c}-n, a, d+1\\2a+1, d\end{array}\!;2\right]
\end{eqnarray*}
upon inversion of the order of summation by hypothesis.

If we now separate the above sum into terms involving even and odd $n$ we obtain 
\[S=\sum_{n=0}^\infty \frac{\phi^{(2n)}(1)}{(2n)!}\,{}_3F_2\left[\!\!\begin{array}{c}-2n, a, d+1\\2a+1, d\end{array}\!;2\right]-\sum_{n=0}^\infty \frac{\phi^{(2n+1)}(1)}{(2n+1)!}\,{}_3F_2\left[\!\!\begin{array}{c}-2n-1, a, d+1\\2a+1, d\end{array}\!;2\right].\]
Finally, using the summations in (\ref{e18a}) and (\ref{e18b}) and noting that $(2n)!=2^{2n}(\fs)_n n!$, $(2n+1)!=2^{2n} (\f{3}{2})_n n!$, we easily arrive at the right-hand side of (\ref{e21}). This completes the proof
of the theorem.\ \ \ \ \ $\Box$

When $d=2a$ it is seen that (\ref{e21}) reduces to Ramanujan's result in (\ref{e11}).

\vspace{0.6cm}

\begin{center}
{\bf 3. \ Examples of Theorem 2}
\end{center}
\setcounter{section}{3}
\setcounter{equation}{0}
\renewcommand{\theequation}{\arabic{section}.\arabic{equation}}
In this section, we provide some examples of different choices for the function $\phi(t)$ appearing in (\ref{e21}).
Throughout this section we let $k$ denote a non-negative integer. 
\bigskip

(a)\ \ First we consider the simplest choice with $\phi(t)=\exp\,(xt)$, where $x$ is an arbitrary variable (independent of $t$). Then $\phi^{(k)}(t)=x^k\,\phi(t)$ which satisfies the conditions for the validity of (\ref{e21}).
Substitution of the derivatives into (\ref{e21}) and identification of the resulting series as hypergeometric functions immediately yields
\bee\label{e31}
e^{-x}{}_2F_2\left[\!\!\begin{array}{c}a, d+1\\2a+1, d\end{array}\!;2x\right]={}_0F_1\left[\!\!\begin{array}{c}-\\a+\fs\end{array}\!;\f{1}{4}x^2\right]
-\frac{(1-2a/d)x}{2a+1}\,{}_0F_1\left[\!\!\begin{array}{c}-\\a+\f{3}{2}\end{array}\!;\f{1}{4}x^2\right].
\ee
which is a result established by a different method in \cite[Theorem 2]{RP}.
In addition, it is interesting to observe that, since
\[{}_0F_1\left[\!\!\begin{array}{c}-\\a+\fs\end{array}\!;\f{1}{4}x^2\right]=\g(a+\fs) (\fs x)^{\fr-a}\,I_{a-\fr}(x),\]
where $I_\nu$ is the modified Bessel function of the first kind, the result (\ref{e31}) can also be written in terms of $I_\nu$.
\bigskip

(b)\ \ If we let $\phi(t)=\cosh\,(xt)$, we have $\phi^{(2k)}(t)=x^{2k} \cosh\,(xt)$ and $\phi^{(2k+1)}(t)=x^{2k+1}\sinh\,(xt)$. Then (\ref{e21}), after a little simplification making use of the identity 
\[(a)_{2k}=(\fs a)_k (\fs a+\fs)_k\, 2^{2k},\]
and letting $d\ra 2d$, reduces to
\bee\label{e32}
{}_3F_4\left[\!\!\begin{array}{c}\vspace{0.1cm}

 \fs a, \fs a+\fs, d+1\\\fs, a+\fs, a+1, d\end{array}\!;x^2\right]=
\cosh x\ {}_0F_1\left[\!\!\begin{array}{c}-\\a+\fs\end{array}\!;\f{1}{4}x^2\right]
-\frac{(1-a/d)}{2a+1}\,x \sinh x\ {}_0F_1\left[\!\!\begin{array}{c}-\\a+\f{3}{2}\end{array}\!;\f{1}{4}x^2\right].
\ee

\bigskip

(c)\ \ If $\phi(t)=(x-t)^{-b}$, where $b$ is an arbitrary parameter and $x>2$, then
\[\phi^{(k)}(t)=\frac{(b)_k}{(x-t)^{b+k}}.\]
From (\ref{e21}), we therefore find
\[\left(\frac{x}{x-1}\right)^{\!-b}\sum_{k=0}^\infty\frac{2^k(a)_k (b)_k(d+1)_k}{(2a+1)_k (d)_kk!}\,x^{-k}=\sum_{k=0}^\infty
\frac{2^{-2k}(b)_{2k}}{(a+\fs)_k k!}\,(x-1)^{-2k}\hspace{4cm}\]
\[\hspace{6cm}-\frac{(1-2a/d)}{2a+1} \sum_{k=0}^\infty \frac{2^{-2k} (b)_{2k+1}}{(a+\f{3}{2})_k k!}\,(x-1)^{-2k-1},\]
which yields
\[\left(\frac{x}{x-1}\right)^{\!-b}{}_3F_2\left[\!\!\begin{array}{c}a, b, d+1\\2a+1, d\end{array}\!;\frac{2}{x}\right]={}_2F_1\left[\!\!
\begin{array}{c}\fs b, \fs b+\fs\\a+\fs\end{array}\!;\frac{1}{(x-1)^2}\right]\hspace{5cm}\]
\[\hspace{6cm}-\frac{(1-2a/d)b}{2a+1} (x-1)^{-1} {}_2F_1\left[\!\!
\begin{array}{c}\fs b+\fs, \fs b+1\\\ \ a+\f{3}{2}\end{array}\!;\frac{1}{(x-1)^2}\right].
\]
If we put $z=x/(1+x)$, this last result becomes 
\[
(1+z)^{-b}{}_3F_2\left[\!\!\begin{array}{c}a, b, d+1\\2a+1, d\end{array}\!;\frac{2z}{1+z}\right]\hspace{7cm}\]
\bee\label{e33}
\hspace{2cm}={}_2F_1\left[\!\!
\begin{array}{c}\fs b, \fs b+\fs\\a+\fs\end{array}\!;z^2\right]
-\frac{(1-2a/d)b}{2a+1}\,z\,{}_2F_1\left[\!\!
\begin{array}{c}\fs b+\fs, \fs b+1\\\ \ a+\f{3}{2}\end{array}\!;z^2\right],
\ee
which has been obtained by different methods in \cite[Theorem 3]{Kim2}.

\bigskip

(d)\ \ Finally, with $\phi(t)=\exp\,(-x^2t^2/4)$, we have \cite[p.~442]{NIST}
\[\phi^{(k)}(t)=(-1)^k(x/2)^k e^{-x^2t^2/4}\,H_k(xt/2),\]
where $H_k$ is the Hermite polynomial of order $k$. Since $H_{2k}(0)=(-1)^k (2k)!/k!$ and $H_{2k+1}(0)=0$, it follows from (\ref{e21}) that
\[e^{x^2/4}\,{}_3F_3\left[\!\!\begin{array}{c}\fs a, \fs a+\fs, d+1\\a+\fs, a+1, d\end{array}\!;-x^2\right]\hspace{7cm}\]
\bee\label{e34}
\hspace{3cm}=\sum_{k=0}^\infty\frac{(x/4)^{2k}}{(a+\fs)_k k!}\,H_{2k}(x/2)+\frac{(1-a/d)}{a+\fs}
\sum_{k=0}^\infty\frac{(x/4)^{2k+1}}{(a+\f{3}{2})_k k!}\,H_{2k+1}(x/2)
\ee
provided $a\neq -\fs, -\f{3}{2}, \ldots \ $, where we have put $d\ra 2d$. 

The above series involving the Hermite polynomials can be expressed in terms of ${}_2F_2$ functions since
\cite[Eqs. (36), (37)]{Kim} 
\[\sum_{k=0}^\infty \frac{(x/4)^{2k}}{(a+\fs)_k k!} \,H_{2k}(x/2)=e^{x^2/4}\,{}_2F_2\left[\!\!\!
\begin{array}{c} \fs a, \fs a+\fs\\a, a+\fs\end{array}\!;-x^2\right],\]
\[\sum_{k=0}^\infty \frac{(x/4)^{2k}}{(a+\f{3}{2})_k k!} \,H_{2k+1}(x/2)=x e^{x^2/4}\,{}_2F_2\left[\!\!\!
\begin{array}{c} \fs a+1, \fs a+\f{3}{2}\\a+\f{3}{2}, a+2\end{array}\!;-x^2\right]\]
to yield
\[{}_3F_3\left[\!\!\begin{array}{c}\fs a, \fs a+\fs, d+1\\a+\fs, a+1, d\end{array}\!;-x^2\right]=
{}_2F_2\left[\!\!\!\begin{array}{c} \fs a, \fs a+\fs\\a, a+\fs\end{array}\!;-x^2\right]
+\frac{x^2(1-a/d)}{4a+2}\, {}_2F_2\left[\!\!\!
\begin{array}{c} \fs a+1, \fs a+\f{3}{2}\\a+\f{3}{2}, a+2\end{array}\!;-x^2\right].\]
We remark that this last result can be derived alternatively by writing $(d+1)_k/(d)_k=1+k/d$ in the series expansion of the ${}_3F_3$ function combined with use of the result for contiguous ${}_2F_2$ functions \cite{Kim}
\[{}_2F_2\left[\!\!\begin{array}{c}\alpha, \beta\\\gamma, \delta\end{array}\!;z\right]
-{}_2F_2\left[\!\!\begin{array}{c}\alpha, \beta\\\gamma, \delta+1\end{array}\!;z\right]
=\frac{\alpha\beta z}{\gamma\delta(\delta+1)}\,
{}_2F_2\left[\!\!\begin{array}{c}\alpha+1, \beta+1\\\gamma+1, \delta+2\end{array}\!;z\right].\]

Finally, the representation (\ref{e34}) may be contrasted with the more general result obtained from Theorem 1 with $\phi(t)=\exp (-x^2t^2/4)$ given in \cite[Eq. (40)]{Kim} (with $x\ra 2x$ and $d\ra 2d$)
\[e^{x^2}\,{}_3F_3\left[\!\!\begin{array}{c}\fs a, \fs a+\fs, d+1\\ \fs b, \fs b+\fs, d\end{array}\!;-x^2\right]
=\sum_{k=0}^\infty \frac{(b-a-1)_k (f+1)_k}{(b)_k (f)_k}\,x^k H_k(x),\]
where $f=2d(b-a-1)/(2d-a)$.
\vspace{0.6cm}

\begin{center}
{\bf 4. \ Extension of the summation (\ref{e19})}
\end{center}
\setcounter{section}{4}
\setcounter{equation}{0}
\renewcommand{\theequation}{\arabic{section}.\arabic{equation}}
We employ the usual convention of writing the finite sequence of parameters $(a_1, \ldots , a_p)$ simply by $(a_p)$ and the product of $p$ Pochhammer symbols by $((a_p))_k\equiv (a_1)_k \ldots (a_p)_k$.
In order to derive our extension of Ramanujan's sum (\ref{e19}) we make use of the generalized Karlsson-Minton summation theorem given below.
\begin{theorem} 
Let $(m_r)$ be a sequence of positive 
integers and $m:= m_1+\cdots +m_r$. The generalized Karlsson-Minton summation theorem is given by \cite{MP4, MS}
\bee\label{e41}
{}_{r+2}F_{r+1}\left[\!\!\begin{array}{c}a, b,
\\c,\end{array}\!\!\!\!\begin{array}{c}(d_r+m_r)\\(d_r)
\end{array}\!;1\right]=
\frac{\g(c) \g(c-a-b)}{\g(c-a) \g(c-b)}\,\sum_{k=0}^m\frac{(-1)^k(a)_k (b)_k C_{k,r}}{(1+a+b-c)_k}
\ee
provided $\Re\,(c-a-b)>m$. 
\end{theorem}
The coefficients $C_{k,r}$ appearing in (\ref{e41}) are defined for $0\leq k\leq m$ by 
\bee\label{e42}
C_{k,r}= \frac{1}{\Lambda}\sum_{j=k}^m \sigma_{j}{\bf S}_j^{(k)},\qquad \Lambda=(d_1)_{m_1}\ldots (d_r)_{m_r},
\ee
with $C_{0,r}=1$, $C_{m,r}=1/\Lambda$. The ${\bf S}_j^{(k)}$ denote the Stirling numbers of the second kind and the $\sigma_j$ $(0\leq j\leq m)$ are generated by the relation
\bee\label{e43}
(d_1+x)_{m_1} \cdots (d_r+x)_{m_r}=\sum_{j=0}^m \sigma_{j}x^j.
\ee 

In \cite{MP3}, an alternative representation for the coefficients $C_{k,r}$ is given as the terminating hypergeometric series of unit argument
\bee\label{e43a}
C_{k,r}=\frac{(-1)^k}{k!}\,{}_{r+1}F_r\left[\!\!\begin{array}{c}-k,\\{}\end{array}
\!\!\!\begin{array}{c}(d_r+m_r)\\(d_r)\end{array}\!;1\right].
\ee
When $r=1$, with $d_1=d$, $m_1=m$, Vandermonde's summation theorem \cite[p.~243]{S} can be used to show that
\bee\label{e44}
C_{k,1}=\left(\!\!\!\begin{array}{c}m\\k\end{array}\!\!\!\right)\,\frac{1}{(d)_k}.
\ee

Our extension of Ramanujan's summation in (\ref{e19}) is given by the following theorem.
\begin{theorem} 
Let $(m_r)$ be a sequence of positive integers and $m:= m_1+\cdots +m_r$. Then, provided
Re$\,(c-a)>m$, we have
\bee\label{e45}
\sum_{k=1}^\infty \frac{(a)_k\,((d_r+m_r))_k}{(c)_k\,((d_r))_k\,k}=\psi(c)-\psi(c-a)+\sum_{k=1}^m
\frac{(-1)^k (a)_k \g(k) C_{k,r}}{(1+a-c)_k},
\ee
where $C_{k,r}$ are the coefficients defined in (\ref{e43}) and (\ref{e43a}).
\end{theorem}
\vspace{0.2cm}

\noindent{\it Proof.}\ \ \ We follow the method of proof given in \cite[p.~25]{B}. If we differentiate logarithmically the left-hand side of (\ref{e41}) with respect to $b$ and then set $b=0$, making use of the simple fact that 
\[\left.\frac{d}{dx} (x)_k\right|_{x=0}=(k-1)!,\qquad k\geq 1,\]
we immediately obtain
\[\sum_{k=1}^\infty \frac{(a)_k\,((d_r+m_r))_k}{(c)_k\,((d_r))_k\,k}\]
when Re$\,(c-a)>m$. Proceeding in a similar manner for the right-hand side of (\ref{e41}), we obtain
\[\psi(c)-\psi(c-a)+\sum_{k=1}^m\frac{(-1)^k (a)_k \g(k) C_{k,r}}{(1+a-c)_k}.\]
This completes the proof of the theorem.\ \ \ \ \ $\Box$
\bigskip

When $r=1$, the coefficients $C_{k,1}$ are given by (\ref{e44}) and we obtain the summation
\bee
\sum_{k=1}^\infty \frac{(a)_k (d+m)_k}{(c)_k (d)_k k}=\psi(c)-\psi(c-a)+m!\sum_{k=1}^m \frac{(-1)^k (a)_k}
{k (m-k)! (1+a-c)_k (d)_k}
\ee
provided Re$\,(c-a)>m$.

\vspace{0.6cm}

\noindent{\bf Acknowledgement:}\ \ \ Y. S. Kim acknowledges the support of the Wonkwang University Research Fund (2014). 

\end{document}